\documentclass[12pt]{article}
\usepackage{personal}
\usepackage{realoracles}

\title{Rationally Querying the Reals}
\jtauthor
\date{April 30, 2023}

\begin{document}\maketitle
\begin{abstract}
A new definition of a real number is that it is a rule which says Yes or No based on whether the real number ought to be in a given rational interval. This is a teaser paper for formalizing, exploring, and generalizing this definition. The full exploration is given in the paper ``Defining Real Numbers as Oracles''. 
\end{abstract}

This is a teaser that gives the basic definitions that other papers explore in detail. These definitions give new ways of completing the rationals into the reals, completing metric spaces, compactifying topological spaces, and a new definition of functions that respects this process. These new objects are called oracles, as they are used to answers Yes or No questions that reveal the item. 

This paper will be brief, stating the definitions and some of the consequences, but not delving much into the reasons why. For that please see the other papers: an overview of the oracles as reals \cite{taylor23over}, the detailed paper on the reals \cite{taylor23main}, topological generalizations \cite{taylor23metric}, elementary analysis using function oracles \cite{taylor23funora}, and some reflections on pedagogical considerations \cite{taylor23edu}. A good resource for alternative definitions of reals can be found in \cite{ittay-2015}. The closest approach to oracles comes from the constructivists using what the oracle papers refer to as fonsis; see \cite{bridger} and \cite{bridges}. 

To know a real number is to be able to answer definitively whether that real number is between any two given rational numbers. With that capability, one can do arbitrarily precise arithmetic with that number. That knowledge is codified in an oracle.

An oracle is a rule, satisfying a few properties below, that takes in an inclusive rational interval and returns a 1 or a 0. The 1 should be returned when an interval ought to contain the real number. This includes the inclusive rational intervals consisting of exactly one point which may be called a singleton. The notation $a:b$ represents the inclusive rational interval. There is no implication of a specific ordering or strict inequality.  Indeed, the  interval $a:a$ represents the singleton consisting of just $a$. The notation $a\lte b$ is used to indicate $a\leq b$ when knowing which one is the lower bound is important. 

In what follows, Yes represents a 1 and No represents a 0 result for the interval. The following are the properties that a rule $R$ must satisfy to be an oracle. 

\begin{enumerate}
    \item Consistency. If an interval contains a Yes interval, then it is a Yes interval. This maintains the illusion of the real number being in the Yes intervals and also leads to this being a unique representative of a real number. 
    \item Existence. There should be a Yes interval. All other properties are conditional so this is the only one that says there is something there. 
    \item Closed. If a rational number is in every Yes interval, then its singleton should also be a Yes interval. This is the other part in making sure that there is only one representative of the real number. It is a primary reason for using inclusive intervals. This also makes arithmetic with rationals have an easy version.
    \item Rooted. There is at most one Yes singleton.  This helps ensure that the oracle is narrowed in on a single number. 
    \item Interval Separation. This is the key property. For a given Yes interval, a point $c$ inside that interval creates two subintervals in addition to its own singleton. The requirement is that exactly one of those intervals is Yes unless $c:c$ is a Yes singleton, in which case all three are Yes since $c$ is in all three.  
\end{enumerate}

If $q$ is contained in every Yes interval, then $q$ is the \textbf{root of the oracle}.  An oracle with a root may be called a rooted oracle or a singleton oracle. Rooted oracles are the rational real numbers. If an oracle is not rooted, then it may be called a neighborly oracle. Neighborly oracles are the irrational real numbers. 

While the first three properties above seem to be basic requirements, for maximality, there is flexibility on the last two. They can be replaced with equivalent properties.

The first equivalent replacement would be the two properties: 1) Two Point Separation, which states that given any two rational numbers in a Yes interval, there exists a Yes interval that does not contain at least one of the two numbers, and 2) Disjointness, which states that if two intervals are disjoint, then at most one of them can be a Yes interval. These two properties are more useful in generalizing to metric spaces, where inclusive balls are considered instead of inclusive intervals. The interval separation does not generalize to that context, but two point separation does. 

The second equivalent replacement would be the two properties: 1) Narrowing, which states that given a rational positive length, there is a Yes interval that is shorter than that length, and 2) Intersection, which is that all Yes intervals intersect. If a set of intervals satisfy these two properties, then they are called a Family of Overlapping, Notionally Shrinking Intervals, or a \textbf{fonsi}. Given a fonsi, there is a unique oracle whose set of Yes intervals contain the fonsi. The oracle is defined via the rule that an interval is considered to be a Yes interval if it contains an intersection of elements of the fonsi. This intersection can be finite or infinite. The infinite intersection is mainly of interest for rooted oracles as the singleton may not be in the fonsi nor included in finite intersections. This is a very practical tool in establishing oracles in a variety of contexts, such as a converging sum. A maximal fonsi is a fonsi in which all intersections of elements in the fonsi that result in rational intervals are in the fonsi and all intervals that contain an element of the fonsi is in the fonsi. A maximal fonsi is exactly the set of Yes intervals for the unique related oracle. An oracle can be viewed as the characteristic function associated with its maximal fonsi. 

Oracles are ordered by looking for disjoint intervals that are Yes intervals for the different oracles. The ordering is then inherited from the disjoint interval ordering. If there are no disjoint Yes intervals, then the oracles are equal. 

Arithmetic is defined by using interval arithmetic of the Yes intervals, with definitive answers coming from the fact that one can choose arbitrarily fine Yes intervals leading to as refined a computation as one would like, in the ideal case.  

Rational real numbers are defined by the rule that an interval is a Yes interval if it contains the rational. Since rationals exist independently of oracles, this makes sense. 

For the positive $n$-th root of a rational $q$, a Yes interval, at least for positive endpoints, is one in which the lower endpoint's $n$-th power is less than or equal to $q$ while the upper endpoint's $n$-th power is greater than or equal to $q$.

More generally, real numbers can often be defined as zeros of a function. The Intermediate Value Theorem setup can often be applied. With appropriate conditions, the oracle here is a rule that yields Yes intervals based on whether the signs of the function values, evaluated at the endpoints, differ. That is, if $a:b$ is the interval of interest, then the oracle says Yes exactly when $f(a)*f(b) \leq 0$.

The bisection algorithm generally used to prove the Intermediate Value can be used on any oracle to produce more narrow Yes intervals. The Separation property guarantees that this works. 

One could also compute mediants instead of midpoints. The mediant algorithm will generate the best rational approximation for a given denominator size as well as producing the continued fraction representation. For rational solutions, the mediant process will produce that rational in a finite number of steps while the midpoint process is not guaranteed to do so. 

The oracles are the real numbers. The least upper bound of a set is identified as the oracle whose Yes intervals are those whose lower endpoint is not above any upper bound of the set while the upper endpoint is an upper bound for the given set. Cauchy sequences have the limit defined as the oracle for which a Yes interval is one which contains the tail of the sequence; singleton limits generally need to have the rational value added to the rule. 

Oracles do generalize to other spaces. The Separation property needs to be replaced with the Two Point Separating property. Intervals are replaced with inclusive balls for completing metric spaces. It is possible to generalize this also to the linear structure topologies of \cite{maudlin}. For Tychonoff spaces, certain closed sets may be used to create a compactification of it. 

Oracles also impact functions. Since an oracle is about Yes intervals, it makes sense to think of functions as not taking in an oracle, but rather as operating on intervals. This leads to the idea of a function oracle being a rule that ought to be saying Yes to a rectangle if the image of the base under the function is contained in the wall of the rectangle. The properties of function oracles are similar to the real oracles, but with the complication that the base must be factored into it. 

By working with the properties, it can be shown that for every oracle $\alpha$ in the base of a Yes rectangle, there is a unique oracle $\beta$ whose Yes intervals are the walls of the Yes rectangles whose base contains $\alpha$. The common notation of $f(\alpha) = \beta$ for this association is therefore justified. 

It can be proven that oracle functions are equivalent to the set of usual real functions that are continuous everywhere except possibly at the rationals. 

Oracles are a new way of perceiving the real numbers. They give a clear and accessible definition while also being able to help with computing out the number. They encourage precision thinking about the imprecision associated with a real number. The concept nicely generalizes not only to other topological spaces, but also to the very notion of what functions on these spaces might mean. They represent how finite beings such as ourselves can meaningfully claim knowledge of an inherently infinite real number.

\medskip

\normalem 
\printbibliography

\end{document}